\documentclass{amsproc}

\usepackage{amssymb}

\usepackage{graphicx}

\usepackage{hyperref}
\usepackage{a4wide}
\usepackage{color}
\usepackage{amsthm}
\usepackage{amsmath}
\usepackage{amscd} \usepackage{verbatim}
\usepackage[all]{xy}
\usepackage{mathtools}
\usepackage{float}
\usepackage{caption}

\newtheorem{theorem}{Theorem}[section]
\newtheorem{corollary}[theorem]{Corollary}

\newtheorem{proposition}[theorem]{Proposition}

\theoremstyle{definition}

\newtheorem{example}[theorem]{Example}

\theoremstyle{remark}
\newtheorem*{remark}{Remark}

\numberwithin{equation}{section}

\newcommand{\la}{\lambda}

\newcommand{\C}{\mathbb{C}}

\newcommand{\F}{\mathbb{F}}

\def\onehalf{\frac{1}{2}}

\def\Leg{\text {\rm L}}
\def\Ftwoone{  _2F_1}
\def\atwoone{ a_{\Leg}}

\def\Etwoone{ E_{\Leg}}

\def\ff{\text {\rm ff}}
\def\cl{\text {\rm cl}}
\def\CL{\text{\rm CL}}

\newcommand{\set}[1]{\left\{ #1 \right\} }

\newcommand{\Aut}{\operatorname{Aut}}
\newcommand{\Mat}{\operatorname{Mat}}
\newcommand{\Spec}{\operatorname{Spec}}

\newcommand{\GL}{\operatorname{GL}}

\DeclarePairedDelimiter{\abs}{\lvert}{\rvert}

\let\set\relax
\DeclarePairedDelimiter{\set}{\{}{\}}

\DeclarePairedDelimiter{\parens}{\lparen}{\rparen}

\newcommand{\Fq}{{\F_q}}
\newcommand{\Zhat}{\hat Z}

\begin{document}

\title[Counting matrix points on certain varieties over finite fields]{Counting matrix points on certain varieties over finite fields}


\author{Yifeng Huang}
\address{Department of Mathematics, University of British Columbia, Vancouver, BC Canada V6T 1Z2}
\curraddr{}
\email{huangyf@math.ubc.ca}
\thanks{Y.H. thanks the AMS--Simons Travel Grant for making this collaboration possible.}

\author{Ken Ono}
\address{Department of Mathematics, University of Virginia, Charlottesville, VA 22904}
\curraddr{}
\email{ken.ono691@virginia.edu}
\thanks{K.O. thanks  the Thomas Jefferson Fund and the NSF (DMS-2002265 and DMS-2055118) for their support.}

\author{Hasan Saad}
\address{Department of Mathematics, University of Virginia, Charlottesville, VA 22904}
\curraddr{}
\email{hs7gy@virginia.edu}
\thanks{}

\subjclass[2020]{33C70; 14Gxx}
\keywords{Hypergeometric functions, Matrix points, Elliptic curves, $K3$ surfaces}

\date{}

\begin{abstract} Classical hypergeometric functions are well-known to play an important role in arithmetic algebraic geometry. These functions offer solutions to ordinary differential equations, and special cases of such solutions are periods of Picard-Fuchs varieties of Calabi-Yau type. Gauss' $_2F_1$  includes the celebrated case of elliptic curves through the theory of elliptic functions. In the 80s, Greene defined finite field hypergeometric functions that can be used to enumerate the number of finite field points on such varieties. We extend some of these results to count finite field ``matrix points." For example, for every $n\geq 1,$ we consider the matrix  elliptic curves
$$
B^2 = A(A-I_n)(A-a I_n),
$$
where $(A,B)$ are commuting $n\times n$ matrices over a finite field $\F_q$ and $a\neq 0,1$ is fixed. Our  formulas are assembled from Greene's hypergeometric functions and $q$-multinomial coefficients. We use these formulas to prove Sato-Tate distributions for the error terms for matrix point counts for these curves and some families of $K3$ surfaces.
\end{abstract}

\maketitle

	\section{Introduction and Statement of Results}

Classical hypergeometric functions are well known to give periods of elliptic  curves. To be precise, if $n$ is a nonnegative integer, then define
$(\gamma)_n$ by
\begin{displaymath}
(\gamma)_n:=\begin{cases} 1 \ \ \ \ \ &{\text {\rm if}}\ n=0,\\
           \gamma(\gamma+1)(\gamma+2)\cdots (\gamma+n-1) \ \ \ \ \
&{\text {\rm if}}\ n\geq 1.
\end{cases}
\end{displaymath}
The {\it classical hypergeometric function} in parameters
$\alpha_1,\dots, \alpha_h,
\beta_1,\dots, \beta_j\in \C$ is defined by
$$_hF^{\cl}_j \left ( \begin{matrix}
\alpha_1& \alpha_2& \dots & \alpha_h\\
\ &\beta_1 &\dots &\beta_j \end{matrix} \ \vrule \ x \right ) :=
\sum_{n=0}^{\infty} \frac{(\alpha_1)_n (\alpha_2)_n
(\alpha_3)_n \cdots
(\alpha_h)_n}{(\beta_1)_n(\beta_2)_n \cdots (\beta_j)_n}\cdot\frac{x^n}{n!}.
$$
Perhaps the most famous example illustrating the role of these functions in geometry involves the Legendre elliptic curves
\begin{equation}
\Etwoone(a):\ \ y^2=x(x-1)(x-a), \ \ \ \ a\in \C\setminus \{0,1\}.
\end{equation}
The theory of elliptic integrals shows, for $0<a<1,$ that
the  function
 $\Ftwoone^{\cl}(x):=\, \Ftwoone^{\cl}\left(\begin{matrix}\onehalf &\onehalf\\ \ &\ 1\
\end{matrix} \ \vrule \ x\right)$ (for example, see
page 184 of \cite{H})  gives the real period $\Omega_{\textrm L}(a)$ of $\Etwoone(a)$ by the formula 
\begin{equation}\label{GaussPeriod}
\Omega_{\textrm L}(a)=\pi \cdot \  \Ftwoone^{\cl}(a).
\end{equation}

There is another kind of hypergeometric function,
the finite field hypergeometric function, that gives further information about these elliptic curves and higher dimensional varieties. These functions count points over finite fields.
To make this precise, we first recall their definition 
which is due to Greene \cite{Greene}.
If $q$ is a prime power and $A$ and $B$ are two Dirichlet characters on $\F_q$
(extended so that $A(0)=B(0)=0$),
then let $\left ( \begin{matrix} A\\ B \end{matrix} \right )$
be the normalized Jacobi sum
\begin{displaymath}
   \left ( \begin{matrix} A \\ B \end{matrix} \right ):=
   \frac{B(-1)}{q} J(A,\overline{B})=
  \frac{B(-1)}{q}\sum_{x\in \Fq} A(x)\overline{B}(1-x).
\end{displaymath}
Here $\overline{B}$ 
is the complex conjugate of $B$.
If $A_{0}, \dots, A_{n},$ and
$B_{1}, \dots, B_{n}$ are characters on $\F_q$,
then the  \textit{finite field
hypergeometric function}
in these parameters is defined by
$$
 _{n+1}F^{\ff}_{n}\left (
\begin{matrix} A_0 & A_1 & \dots & A_n \\ \ &  B_1 & \dots & B_n \end{matrix}
\ \vrule \  x \right )_q :=\frac{q}{q-1}
  \sum_{\chi} \left ( \begin{matrix} A_0 \chi \\ \chi \end{matrix} \right )
              \left ( \begin{matrix} A_1 \chi \\ B_1 \chi \end{matrix} \right )
              \cdot \cdot \cdot
              \left ( \begin{matrix} A_n \chi \\ B_n \chi \end{matrix} \right )
              \chi(x).
$$
Here $\sum_{\chi}$ denotes
the  sum over all characters $\chi$ of $\F_q$.

It has been observed by many authors (see  \cite{Greene}, \cite{GreeneStanton}, \cite{ISS}, \cite{Koike},  \cite{OnoTAMS}, and \cite{Rouse}, to name a few) that the Gaussian analog of a classical hypergeometric series with rational parameters is obtained by replacing each $\frac{1}{n}$ with a character $\chi_n$ of order $n$ (and $\frac{a}{n}$ with $\chi_n^a$). Let $q$ be a power of an odd prime, $\epsilon_q$ be the trivial character on $\F_q$
and let $\phi_q$ be the character of order $2$. Then the  finite field analog
of $\Ftwoone^{\cl}(x)$ is 
\[
\Ftwoone^{\ff}(x)_q:= {_2F^{\ff}_1}\left (
\begin{matrix} \phi_q &  \phi_q  \\ \ &
\epsilon_q   \end{matrix}
\ \vrule \ x \right )_q.
\]
More generally, we let
\begin{equation}
 _{n+1}F^{\ff}_{n}(x)_q:=
\ _{n+1}F^{\ff}_{n}\left (
\begin{matrix} \phi_q &  \phi_q & \dots \phi_q \\ \ &
\epsilon_q & \dots \epsilon_q \end{matrix}
\ \vrule \ x \right )_q.
\end{equation}
M. Koike proved \cite{Koike} that if  $p$ is an odd prime, $q$ is a power of $p$ and $a\in\F_q\setminus\{0,1\},$  then
\begin{equation}\label{Koiketrace}
    \Ftwoone^{\ff}(a)_q=-\frac{\phi_q(-1)}{q}\cdot \atwoone(a;q),
\end{equation}
where $q+1-\atwoone(a;q)$ counts the number of $\F_q$-points  on $\Etwoone(a)$. This expression is the finite field analogue of Gauss' period formula (\ref{GaussPeriod}).

Motivated by (\ref{GaussPeriod}) and (\ref{Koiketrace}), it is natural to ask whether other finite field hypergeometric function evaluations give point counts for other varieties. This is indeed the case, and perhaps the most beautiful example involves  the analog of the celebrated classical Clausen identity  \cite{Clausen}
\begin{equation}
{ _3F^{\cl}_2\left(\begin{matrix}
	b+c & 2b & 2c \\
	\ & b+c+\frac{1}{2} & 2b+2c
\end{matrix} \ \vrule \ x \right)={_2F^{\cl}_1}\left(\begin{matrix}
b & c \\ 
\ & b+c+\frac{1}{2}
\end{matrix} \ \vrule \ x \right)^2.}
\end{equation}
Using this identity, D. McCarthy \cite{McCarthy} proved that if $a>0,$ then 
\begin{equation*}
{_3F_2^\cl}\left(\begin{matrix}
	\frac{1}{2} & \frac{1}{2} & \frac{1}{2} \\
	\ & 1 & 1 
\end{matrix} \ \vrule \ \frac{a}{a+1}\right)=\frac{\sqrt{1+a}}{\pi^2}\cdot \Omega_{\CL}(a)^2,
\end{equation*}
where $\Omega_{\CL}(a)$ is the real period of the Clausen elliptic curve
$$
E_{\CL}(a):\ \ y^2=(x-1)(x^2+a).
$$

In the finite field case, the second author proved [Theorem 5 of \cite{OnoTAMS}] that if $\F_q$ is a finite field of characteristic $\text{char}(\F_q)\geq 3$ and $a\in\F_q\setminus\{0,-1\},$ then
\begin{equation}\label{Clausen3F2}
q+q^2\phi_q(a+1)\cdot{_3F^{\ff}_2}\left(\frac{a}{a+1}\right)_q=a_{\CL}(a;q)^2=q^2\cdot {_2F^{\ff}_1}\left(\frac{1-\sqrt{-a}}{1+\sqrt{-a}}\right)_q^2,
\end{equation}
where $q+1-a_{\CL}(a;q)$ is the number of $\F_q$ points on $E_{\CL}(a),$ and where the second equality holds whenever $-a$ is a square in $\F_q$ by noting that $E_{\CL}(a)$ is a quadratic twist of $E_{\Leg}\left(\frac{1-\sqrt{-a}}{1+\sqrt{-a}}\right)$ and using (\ref{Koiketrace}). This equality is an analogue of a special case of Clausen's identity.
 Furthermore, this identity can be interpreted in terms of $K3$ surfaces whose function fields are given by
$$
X_a:\ \ s^2=xy(x+1)(y+1)(x+a y),
$$
where $a\in\F_q\setminus\{0,-1\}.$ In this notation, it is known (see Theorem 11.18 of \cite{onoWebOfModularity} and Proposition 4.1 of \cite{ahlgrenOnoPenniston}) that
\begin{equation}\label{K3FqCase}
|X_a(\F_q)|=1+q^2+19q+q^2\cdot{_3F_2^{\ff}}(-a)_q.
\end{equation}

In this note we show that the hypergeometric identities (\ref{Koiketrace}) and (\ref{K3FqCase}),  combined with the combinatorial input from partitions and $q$-multinomial coefficients, 
  count suitable  ``matrix points'' on  these curves and surfaces. To make this precise, we first introduce some notation. If $n,m$ are positive integers and $K$ is a field, then let $C_{n,m}(K)$ denote the set of pairwise-commuting $m$-tuples of $n\times n$-matrices over $K.$ Due to the noncommutativity of matrix multiplication, geometric problems related to studying matrix rational points on curves and higher varieties only make sense when the matrices are commuting, that is, when the matrix points are in $C_{n,m}(K).$ We will be interested in counting tuples in $C_{n,m}(K)$ which satisfy the equations defining some affine varieties. More precisely, we will consider the sets
	$$
	\{(A,B)\in C_{n,2}(\F_q) : B^2 = A(A-I_n)(A-a I_n)\}
	$$
	and
	$$
	\{(A,B,C)\in C_{n,3}(\F_q) : C^2 = AB(A+I_n)(B+I_n)(A+a B), C\in\GL_n(\F_q)\}
	$$
	as matrix analogues of (the smooth affine parts of) the Legendre elliptic curves $E_{\Leg}$ and the K3 surfaces $X_a$ considered above. 
	
	To express our results, we introduce some notation.
	 If $\lambda$ is a partition of a nonnegative  integer $k,$ we write $n(\lambda; i)$ to denote the number of times $i$ is repeated in $\lambda.$ Furthermore, we write $|\lambda|=k,$ and write $l(\lambda)=\sum n(\lambda;i)$ to denote the number of parts of $\lambda.$ Additionally, we introduce certain polynomials in $q.$ More precisely, if $z$ and $q$ are any complex numbers, and $n$ is any positive integer, then we define the $q$-Pochhammer symbol
	 \begin{equation}
	 	(z;q)_n:=(1-z)(1-zq)\dots(1-zq^{n-1})
	 \end{equation}
 	with $(z;q)_0=1.$ 
 	The series expansion of $(q;q)_5$ will play a special role in our results. For clarity of results, we define $b_r$ for integers $r\geq 0$ by
 	\begin{equation}\label{BRCoefficients}
 	\sum\limits_{r=0}^\infty b_rq^r=(q;q)_\infty^5=\parens*{1+\sum_{m=1}^\infty (-1)^m q^{m(3m-1)/2}(1+q^m)}^5,
 	\end{equation}
 	where the last expression due to Euler's pentagonal number theorem \cite[Corollary 1.7]{andrews} allows explicit computation for each $b_r$.
 	
	 Finally, for an integer $n\geq 0$ and $m_1+\ldots+m_k=n$ a partition of $n,$  we define the $q$-multinomial factor
	$$
	\binom{n}{m_1,m_2,\ldots,m_k}_q:=\frac{(q;q)_n}{(q;q)_{m_1}(q;q)_{m_2}\ldots(q;q)_{m_k}}.
	$$
	It is known that $\binom{n}{m_1,\ldots,m_k}_q$ is a monic polynomial in $q$ and that $\binom{n}{m_1,\ldots,m_k}_q$ approaches the usual multinomial coefficient as $q\to 1.$
	
	We start by expressing the number of commuting matrices on a Legendre elliptic curve. More precisely, if $n$ is a positive integer, $q$ is a prime power and $a\in\F_q,$ we let
	\begin{equation}
	N_{n,2}(a;q) := |\{(A,B)\in C_{n,2}(\F_q) : B^2 = A(A-I_n)(A-a I_n)\}|.
	\end{equation}
	In this notation, we have the following theorem that determines these counts, and also explains the connection with the classical $_2F_1^{\cl}$-hypergeometric function.

\begin{theorem}\label{EllipticCurveCase}
If $q=p^r$ is a prime power with $p\geq 3$ and  $a\in\F_q\setminus\{0,1\},$ then
$$
N_{n,2}(a; q)=P(n,0)_q-\sum\limits_{k=1}^{n}\phi_{q^k}(-1)\cdot P(n,k)_q\cdot {_2F^{\ff}_1}(a)_{q^k},
$$
where
$$
P(n,k)_q := (-1)^kq^{n(n-k)+\frac{k(k+1)}{2}}\sum\limits_{s=0}^{\lfloor\frac{n-k}{2}\rfloor} q^{2s(s-n+k)}\binom{n}{s,n-k-2s,k+s}_q.
$$
Moreover, $P(n,k)_q$ is a polynomial in $q$ with leading term $(-1)^k\cdot  q^{n^2-\frac{k(k-1)}{2}}$ and
$$
\lim\limits_{q\to 1} P(n,k)_q = (-1)^k \binom{n}{k}\cdot {_2F^{\cl}_1}\left(\begin{matrix}
	\frac{k-n}{2} & \frac{k+1-n}{2}\\
	\ & k+1
\end{matrix}\ \vrule \ 4\right).
$$
\end{theorem}
\begin{remark}
Although the number theoretic results in this paper require that each $q$ is a prime power, we included the limiting behavior
$\lim_{q\rightarrow 1} P(n,k)_q$ to illustrate that these number theoretic quantities are $q$-analogues of a classical $ _2F_1^{\cl}$ evaluation.
\end{remark}

As a corollary, we consider the matrix analog of the Sato--Tate distribution for point counts for elliptic curves over finite fields. In direct analogy, we find that the limiting distribution of the ``random part'' of matrix point counts on Legendre elliptic curves is semicircular. More precisely, if $n$ is a positive integer and $q$ is a prime power, then we let 
\begin{equation}
a_{\Leg,n}(a; q):=N_{n,2}(a; q)-P(n,0)_{q}.
\end{equation}
In this notation, we have the following result.
\begin{corollary}\label{LegendreDistributionMatrices}
If $-2\leq b<c\leq 2$ and $n$ and $r$ are fixed positive integers, then we have
$$
\lim\limits_{p\to\infty}\frac{|\{a\in\F_{p^r} : p^{\frac{r}{2}-rn^2}a_{\Leg,n}(a; p^r)\in [b,c]\}|}{p^r} = \frac{1}{2\pi}\int_b^c \sqrt{4-t^2}dt.
$$
\end{corollary}

\begin{example}
	For the prime $p=93283,$ we compare the histogram of the distribution of $p^{-7/2} a_{L,2}(a; p)$ for $a\in\F_{p}$ with the limiting distribution.
	\begin{center}
		\begin{table}[H]
		 \includegraphics[height=40mm]{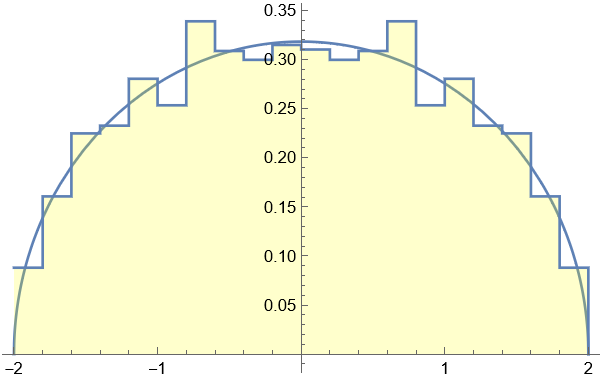}
			\caption*{ $p^{-7/2}a_{L,2}(a;p)$ histogram for $p=93283$}
		\end{table}
	\end{center}
\end{example}

We also consider the matrix version of the $K3$ surfaces described above.
If $n$ is a positive integer, $q$ is a prime power, and $a\in\F_q,$ then we let 
\begin{equation}\label{eq:matrix_point_k3}
N_{n,3}(a;q) := |\{(A,B,C)\in C_{n,3}(\F_q) : C^2 = AB(A+I_n)(B+I_n)(A+a B), C\in\GL_n(\F_q)\}|.
\end{equation}
In this notation, we have the following theorem that gives matrix point counts in terms of the $_3F_2^{\ff}$-hypergeometric function, 6-tuples of integer partitions, the coefficients $b_r$ in (\ref{BRCoefficients}), and $q$-multinomial coefficients.

\begin{theorem}\label{K3Case}
If $q=p^r$ is a prime power with $p\geq 3$ and $a\in\F_q\setminus\{0,-1\},$ then we have
$$
N_{n,3}(a; q)=R(n,\phi_q(a+1))_q+Q\left(n,0,\phi_q(a+1)\right)_q+\sum\limits_{k=1}^n Q\left(n,k,\phi_q(a+1)\right)_q\cdot {_3F^{\ff}_2}\left(-a\right)_{q^k},
$$
where
\begin{align*}
Q(n,k,\gamma)_q :=\ & q^{\frac{n(n-1)}{2}}\sum\limits_{r=k}^n b_{n-r}\sum\limits_{\substack{\lambda_1,\ldots,\la_6 \\ |\lambda_1|+\ldots+|\lambda_6|=r \\ l(\lambda_5)-l(\lambda_6)=k}}  q^{2l(\lambda_3)+l(\la_4)+2l(\la_5)}\gamma^{l(\lambda_4)}(-1)^{n-m(\lambda_1,\ldots,\lambda_6)}\\ 
& (q;q)_{n-m(\lambda_1,\ldots,\lambda_6)}\cdot q^{\sum\frac{n(\lambda_i,j)(n(\lambda_i,j)-1)}{2}}\cdot\binom{n}{n(\lambda_i,j), n-m(\lambda_1,\ldots,\lambda_6)}_q
\end{align*}
and
\begin{align*}
R(n,\gamma)_q :=\ & -q^{\frac{n(n-1)}{2}}\sum\limits_{k=1}^n q^k\sum\limits_{r=1}^n b_{n-r}\sum\limits_{\substack{\lambda_1,\ldots,\la_6 \\ |\lambda_1|+\ldots+|\lambda_6|=r \\ l(\lambda_5)-l(\lambda_6)=k}} \gamma^k q^{2l(\lambda_3)+l(\la_4)+2l(\la_6)}\gamma^{l(\lambda_4)}(-1)^{n-m(\lambda_1,\ldots,\lambda_6)}\\ 
& (q;q)_{n-m(\lambda_1,\ldots,\lambda_6)}\cdot q^{\sum\frac{n(\lambda_i,j)(n(\lambda_i,j)-1)}{2}}\cdot\binom{n}{n(\lambda_i,j), n-m(\lambda_1,\ldots,\lambda_6)}_q
\end{align*}
with $\lambda_1,\ldots,\lambda_6$ being partitions and $m(\lambda_1,\ldots,\lambda_6)=\sum\limits_{i=1}^6 l(\lambda_i).$
Moreover, $Q(n,k,\gamma)_q$ is a polynomial in $q$ with leading term $q^{n^2+n}$ and 
\begin{equation}\label{crazy}
\lim\limits_{q\to 1} Q(n,k,\gamma)_q = \binom{n}{k} \cdot \gamma^{n-k}\sum\limits_{s=0}^{n-k} \binom{n-k}{s}\cdot\left(\frac{3}{\gamma}\right)^s{_2F^{\cl}_1}\left(\begin{matrix}
	\frac{k-n+s}{2} & \frac{k+1-n+s}{2}\\
	\ & k+1
\end{matrix}\ \vrule \ \frac{4}{\gamma^2}\right),
\end{equation}
when $\gamma\neq 0,$ and 
$$
\lim\limits_{q\to 1}Q(n,k,0)_q=\binom{n}{k}\cdot 3^{n-k}\cdot {_2F^{\cl}_1}\left(\begin{matrix}
	\frac{k-n}{2} & \frac{k+1-n}{2}\\
	\ & k+1
\end{matrix}\ \vrule \ \frac{4}{9}\right).
$$
\end{theorem}

\begin{remark} Although the number theoretic results in this paper require that each $q$ is a  prime power, we included the limiting behavior $\lim_{q\rightarrow 1} Q(n,k,\gamma)_q$
to illustrate that these formulas are $q$-analogues of natural weighted sums of $ _2F_1^{\cl}$ evaluations (i.e. see (\ref{crazy}))
that enjoy nice recurrence relations. Namely, if we let
$$
F(z;m,k):=z^m\sum\limits_{s=0}^m \binom{m}{s}\left(\frac{3}{z}\right)^s{_2F^{\cl}_1}\left(\begin{matrix}
	\frac{s-m}{2} & \frac{s-m+1}{2}\\
	\ & k+1
\end{matrix}\ \vrule \ \frac{4}{z^2}\right),
$$
then we have 
$$
z\frac{d}{dz}F(z;m,k)=mF(z;m,k)-3mF(z;m-1,k)-\frac{2m(m-1)}{k+1}F(z;m-2,k+1),
$$
where $F(z;0,k)=1$ for $k\geq 0.$ This follows immediately from
the fact that
	$$
	\frac{d}{dz}{_2F^{\cl}_1}\left(\begin{matrix}
		a & b \\ 
		\ & c
	\end{matrix} \ \vrule \ z \right)=\frac{ab}{c}\cdot {_2F^{\cl}_1}\left(\begin{matrix}
	a+1 & b+1 \\ 
	\ & c+1
	\end{matrix} \ \vrule \ z \right).
	$$
\end{remark}

Theorem~\ref{K3Case} allows us to determine the Sato--Tate type limiting distribution of the ``random part'' of matrix point counts on the $K3$ surfaces $X_a.$ More precisely, if $n$ is a positive integer and $q$ is a prime power, then we let
\begin{equation}
A_n(a; q) := N_{n,3}(a;q)-Q(n,0,\phi_{q}(a+1))_{q}-R(n,\phi_{q}(a+1))_{q}.
\end{equation}
In this notation, we have the following result.
\begin{corollary}\label{K3DistributionMatrices}
 	If $-3\leq b<c\leq 3$ and $n$ and $r$ are fixed positive integers, then we have 
	$$
	\lim\limits_{p\to\infty}\frac{\{a\in\F_{p^r} : p^{r-rn^2-rn} A_n(a ; p^r)\in [b,c]\}}{p^r} = \frac{1}{4\pi}\int_b^c f(t)dt,
	$$
	where
	\begin{equation*}
		f(t):=\begin{cases}
			\sqrt{\frac{3-|t|}{1+|t|}} \ \ \ \ &\ \ \ \ {\text {\it if}}\  1<|t|<3, \\ \\
			\sqrt{\frac{3-t}{1+t}}+\sqrt{\frac{3+t}{1-t}} &\ \ \ \  {\text{{\it if}}}\  |t|<1,\\ \\
			0 &\ \ \ \ \text{otherwise}.
		\end{cases}
	\end{equation*}
\end{corollary}

\begin{example}
	For the prime $p=93283,$ we compare the histogram of the distribution of $p^{-5} A_{2}(a; p)$ for $a\in\F_{p}$ with the limiting distribution.
	\begin{center}
		\begin{table}[H]
 \includegraphics[height=40mm]{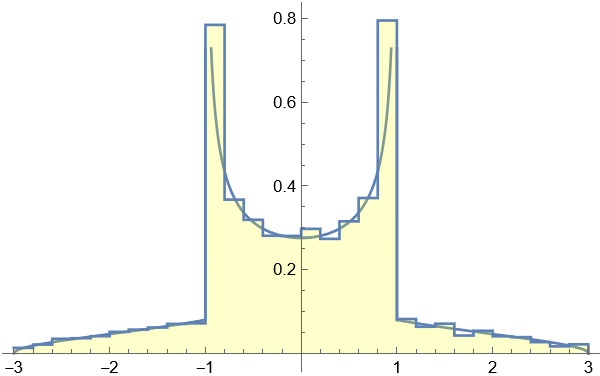}
			\caption*{ $p^{-5}A_{2}(a;p)$ histogram for $p=93283$}
		\end{table}
	\end{center}
\end{example}

\begin{remark}
	The results of this paper can be extended almost {\it mutatis mutandis} to other hypergeometric families of varieties of dimensions $1$ and $2$ such as those introduced by Beukers, Cohen, and Mellit in \cite{BCM}.
\end{remark}

This paper is organized as follows. In Section 2 we recall properties of zeta functions for curves and surfaces in the commuting matrix situation. These results \cite{huang} are due to the first author. In Section 3 we recall results of the second two authors, which we then combine with these zeta functions to obtain our results.

\section*{Acknowledgements} \noindent
	The authors thank the anonymous referees for their remarks and corrections.
		
	\section{Some zeta functions}

	Let $q$ be a prime power. Recall that $\GL_n(\Fq)$ is the group of $n\times n$ invertible matrices over the finite field $\Fq$ with $q$ elements. 
	It will be repetitively used in this paper that
	\begin{equation}\label{GLnCount}
	\abs{\GL_n(\Fq)}= (-1)^n q^{\frac{n(n-1)}{2}} (q;q)_n.
	\end{equation}
	
	Now, let $X=\Spec R$ be an affine variety over $\Fq$. Say
	\begin{equation}
 	R:=\frac{\Fq[T_1,\dots,T_m]}{(f_1,\dots,f_r)}.
 	\end{equation} 	
	
	Following the work \cite{huang} of the first author, we define the set of $n\times n$ \emph{matrix points} on $X$ as the set of commuting tuples of matrices satisfying the defining equations for $X$:
	\begin{equation}\label{eq:matrix_point_def}
	C_n(X):=\set[\bigg]{\underline{A}=(A_1,\dots,A_m)\in \Mat_n(\Fq)^m: [A_i,A_j]=0, f_i(\underline{A})=0}.
	\end{equation}
	
	Note that $C_1(X)\cong X(\Fq)$. Though not needed in this paper, it is worth pointing out that the cardinality of $C_n(X)$ is independent of the choice of defining equations for $X$; in fact, by comparing \cite[Eq.\ 4.2]{huang} and \cite[Eq.\ 4.15]{huang}, there is an equation-free equivalent characterization for the cardinality of $C_n(X)$:
	\begin{equation}
	\frac{\abs{C_n(X)}}{\abs{\GL_n(\Fq)}} = \sum_{\dim_{\Fq} H^0(X;M)=n} \frac{1}{\abs{\Aut M}},
	\end{equation}
	where the sum ranges over all isomorphism classes of zero-dimensional coherent sheaves on $X$ of degree $n$. This characterization also makes $\abs{C_n(X)}$ well-defined for any variety $X$ over $\Fq$.
	
	The number of matrix points on a smooth curve or a smooth surface is given by infinite product formulas for a zeta function associated to it. For any (affine) variety $X$ over $\Fq$, consider its \emph{Cohen--Lenstra series} (terminology of \cite{huang}):
	\begin{equation}
	\Zhat_X(t):=\sum_{n=0}^\infty \frac{\abs{C_n(X)}}{\abs{\GL_n(\Fq)}}t^n,
	\end{equation}
	and recall the local zeta function
	\begin{equation}
	Z_X(t):=\exp\parens*{\sum_{n=1}^\infty \frac{\abs{X(\F_{q^n})}}{n}\, t^n}.
	\end{equation}

	\begin{proposition}[{\cite[Proposition 4.6(a)]{huang}}]\label{prop:curve_zeta_formula}
	If $X$ is a smooth curve over $\Fq$, then
	\begin{equation}
	\Zhat_X(t)=\prod_{j\geq 1} Z_X(tq^{-j}).
	\end{equation}
	\end{proposition}
	
	\begin{proposition}[{\cite[Proposition 4.6(b)]{huang}}]\label{prop:surface_zeta_formula}
	If $X$ is a smooth surface over $\Fq$, then
	\begin{equation}
	\Zhat_X(t)=\prod_{i,j\geq 1} Z_X(t^i q^{-j}).
	\end{equation}
	\end{proposition}
	
	Proposition \ref{prop:curve_zeta_formula} is essentially due to Cohen and Lenstra \cite{CohenLenstra}, and Proposition \ref{prop:surface_zeta_formula} is essentially due to the Feit--Fine formula \cite{FeitFine} for counting commuting matrices and ideas of Bryan and Morrison \cite{BryanMorrison}. We remark that  both formulas heavily exploit the local geometry of $X$, namely, smoothness of dimension $1$ or $2$. In fact, in light of the main theorem of \cite{huang}, Proposition \ref{prop:curve_zeta_formula} ceases to hold if $X$ is a multiplicative reduction of an elliptic curve over a number field (but holds if it is a good reduction).

	\section{Proofs of Theorems~\ref{EllipticCurveCase} and \ref{K3Case}}
	Here we use the results of the previous section to prove Theorems~\ref{EllipticCurveCase} and \ref{K3Case}  and their corollaries.
	
	\subsection{Proof of Theorem~\ref{EllipticCurveCase}}
	
	 Fix a prime power $q=p^r$ with $p\geq 3$ and $r\geq 1,$ and fix $a\in\F_q\setminus\{0,1\}.$ Then, denoting by $X$ the affine part of $E_{\Leg}(a),$ Theorem V.2.4 of \cite{silvermanEllipticCurves} states that
	$$
	Z_X(t) = \frac{(1-\alpha t)(1-\overline{\alpha}t)}{1-qt},
	$$
	where $\alpha$ and $\overline{\alpha}$ are the eigenvalues of Frobenius acting on the Tate module of $E_{\Leg}(a).$ Note that there is a missing factor of $\frac{1}{1-t}$ in this expression since we are only considering the affine part of $E_{\Leg}(a).$ 
	
	By Proposition~\ref{prop:curve_zeta_formula}, we then have that 
	$$
	\hat{Z}_X(t) = \prod\limits_{j\geq 1} \frac{(1-\alpha tq^{-j})(1-\overline{\alpha}tq^{-j})}{1-tq^{1-j}}.
	$$
	It is well-known due to Euler \cite[Corollary 2.2]{andrews} that  
	\begin{equation}
	\prod\limits_{j\geq 1}(1-ctq^{-j}) = \sum\limits_{m\geq 0} \frac{(ct)^m}{(q;q)_m}
	\end{equation}
	and
	\begin{equation}\label{Euler2}
		\prod\limits_{j\geq 1}(1-ctq^{-j})^{-1}=\sum\limits_{m\geq 0} \frac{(-1)^m q^{m(m-1)/2}\cdot (ct)^m}{(q;q)_m}.
	\end{equation}
	This implies that 
	$$
	\hat{Z}_X(t)=\left(\sum\limits_{r\geq 0}\frac{(\alpha t)^r}{(q;q)_r}\right)\cdot \left(\sum\limits_{s\geq 0}\frac{(\overline{\alpha} t)^s}{(q;q)_s}\right)\cdot\left(\sum\limits_{u\geq 0} \frac{(-1)^u q^{u(u+1)/2}\cdot t^u}{(q;q)_u}\right).
	$$
	By the definition of $\hat{Z}_X(t)$ and by ~(\ref{GLnCount}), we then have 
	$$
	N_{n,2}(a;q) = (-1)^n q^{n(n-1)/2} (q;q)_n\cdot \sum\limits_{\substack{r+s+u=n \\ r,s,u\geq 0}} \frac{\alpha^r\overline{\alpha}^s (-1)^u q^{\frac{u(u+1)}{2}}}{(q;q)_r(q;q)_s(q;q)_u}.
	$$
	Furthermore, again by Theorem V.2.4 of  \cite{silvermanEllipticCurves}, we have that $\alpha\overline{\alpha}=q$ and therefore, we can rewrite this sum as 
	$$
	N_{n,2}(a;q) = (-1)^nq^{n(n-1)/2}(q;q)_n\sum\limits_{\substack{r+s+u=n \\ r,s,u\geq 0}} \frac{(-1)^u\alpha^{r-s}q^{s+\frac{u(u+1)}{2}}}{(q;q)_r(q;q)_s(q;q)_u}.
	$$
	Dividing this sum according to the value of $r-s,$ we then have
	\begin{align*}
		&N_{n,2}(a;q)=(-1)^nq^{n(n-1)/2}(q;q)_n\sum\limits_{k}\alpha^k\cdot \sum\limits_{\substack{s\geq 0\\ r=s+k\geq 0\\ u=n-2s-k\geq 0}}(-1)^uq^{s+u(u+1)/2}\frac{1}{(q;q)_r(q;q)_s(q;q)_u} \\
		&=(-1)^nq^{n(n-1)/2}(q;q)_n\sum\limits_k\alpha^k\sum\limits_{s=\max\{-k,0\}}^{\lfloor\frac{n-k}{2}\rfloor}(-1)^{n-k}q^{s+\frac{(n-2s-k)(n-2s-k+1)}{2}}\cdot\frac{1}{(q;q)_s(q;q)_{s+k}(q;q)_{n-2s-k}} \\
		&=q^{n(n-1)/2}\sum\limits_k \alpha^k(-1)^k\sum\limits_{s=\max\{-k,0\}}^{\lfloor\frac{n-k}{2}\rfloor}q^{\frac{k^2}{2}-kn-\frac{k}{2}+\frac{n^2}{2}+\frac{n}{2}}\cdot q^{2ks-2ns+2s^2}\cdot\frac{(q;q)_n}{(q;q)_s(q;q)_{s+k}(q;q)_{n-2s-k}} \\
		&=q^{n(n-1)/2}\cdot\sum\limits_{s=0}^{\lfloor\frac{n}{2}\rfloor} q^{\frac{n^2}{2}+\frac{n}{2}+2s^2-2ns}\cdot\frac{(q;q)_n}{(q;q)_s(q;q)_s(q;q)_{n-2s}} \\
		&+q^{n(n-1)/2}\sum\limits_{k>0}\alpha^k(-1)^k\sum\limits_{s=0}^{\lfloor\frac{n-k}{2}\rfloor}q^{\frac{k^2}{2}-kn-\frac{k}{2}+\frac{n^2}{2}+\frac{n}{2}}\cdot q^{2ks-2ns+2s^2}\cdot\frac{(q;q)_n}{(q;q)_s(q;q)_{s+k}(q;q)_{n-2s-k}}  \\
		&+q^{n(n-1)/2}\sum\limits_{k<0}q^k\overline{\alpha}^{-k}(-1)^k\sum\limits_{s=-k}^{\lfloor\frac{n-k}{2}\rfloor}q^{\frac{k^2}{2}-kn-\frac{k}{2}+\frac{n^2}{2}+\frac{n}{2}}\cdot q^{2ks-2ns+2s^2}\cdot\frac{(q;q)_n}{(q;q)_s(q;q)_{s+k}(q;q)_{n-2s-k}}.
	\end{align*}
	Replacing $k$ in the last sum with $-k$ and then $s$ with $s+k,$ we then have
	\begin{align*}
		N_{n,2}(a;q)&=P(n,0)_q+\sum\limits_{k>0}\alpha^k(-1)^k\sum\limits_{s=0}^{\lfloor\frac{n-k}{2}\rfloor}q^{\frac{k^2}{2}-kn-\frac{k}{2}+n^2}\cdot q^{2ks-2ns+2s^2}\cdot\frac{(q;q)_n}{(q;q)_s(q;q)_{s+k}(q;q)_{n-2s-k}}  \\
		&+\sum\limits_{k>0}\overline{\alpha}^k(-1)^k\sum\limits_{s=k}^{\lfloor\frac{n+k}{2}\rfloor} q^{\frac{k^2}{2}+kn-\frac{k}{2}+n^2}q^{-2ks-2ns+2s^2}\cdot\frac{(q;q)_n}{(q;q)_s(q;q)_{s-k}(q;q)_{n-2s+k}} \\
		&=P(n,0)_q+\sum\limits_{k>0}\alpha^k(-1)^k\sum\limits_{s=0}^{\lfloor\frac{n-k}{2}\rfloor}q^{\frac{k^2}{2}-kn-\frac{k}{2}+n^2}\cdot q^{2ks-2ns+2s^2}\cdot\frac{(q;q)_n}{(q;q)_s(q;q)_{s+k}(q;q)_{n-2s-k}}  \\
		&+\sum\limits_{k>0}\overline{\alpha}^k(-1)^k\sum\limits_{s=0}^{\lfloor\frac{n-k}{2}\rfloor}q^{\frac{k^2}{2}+kn-\frac{k}{2}+n^2}q^{-2nk+2ks-2ns+2s^2}\cdot\frac{(q;q)_n}{(q;q)_s(q;q)_{s+k}(q;q)_{n-2s-k}} \\
		&=P(n,0)_q+\sum (-1)^k(\alpha^k+\overline{\alpha}^k)\cdot  q^{\frac{k^2}{2}-kn-\frac{k}{2}+n^2}\sum\limits_{s=0}^{\lfloor\frac{n-k}{2}\rfloor}q^{2ks-2ns+2s^2}\cdot\frac{(q;q)_n}{(q;q)_s(q;q)_{s+k}(q;q)_{n-2s-k}}.
	\end{align*}
	Since $a_{\Leg}(a;q^k) = \alpha^k + \overline{\alpha}^k,$ (\ref{Koiketrace}) implies that  
	$$
	N_{n,2}(a;q) = P(n,0)_q-\sum\limits_{k=1}^n \phi_{q^k}(-1)\cdot P(n,k)_q\cdot {_2F^{\ff}_1}(a)_{q^k},
	$$
	where
	$$
	P(n,k)_q=(-1)^k\cdot q^{n(n-k)+\frac{k(k+1)}{2}}\sum\limits_{s=0}^{\lfloor\frac{n-k}{2}\rfloor}q^{2s(s-n+k)}\cdot\frac{(q;q)_n}{(q;q)_s(q;q)_{s+k}(q;q)_{n-2s-k}}.
	$$
	
	The leading coefficient of $P(n,k)_q$ is clear from the expression. Since the $q$-multinomial approaches the usual multinomial as $q\to 1,$ we have
	$$
	\lim\limits_{q\to 1} P(n,k)_q = (-1)^k\sum\limits_{s=0}^{\lfloor\frac{n-k}{2}\rfloor} \binom{n}{s,n-k-2s,k+s} = (-1)^k\binom{n}{k}\sum\limits_{s=0}^{\lfloor\frac{n-k}{2}\rfloor}\frac{k!(n-k)!}{s!(k+s)!(n-k-2s)!}.
	$$
	It is easy to see by induction that if $m$ and $s$ are integers with $m<0$ and $2s+m\leq 0,$ then
	$$
	\left(\frac{m}{2}\right)_s\left(\frac{m}{2}+\frac{1}{2}\right)_s = \frac{(-m)!}{(-m-2s)!4^s}.
	$$
	Furthermore, it is evident by definition that for $k\geq 0,$ we have $(k+1)_s=\frac{(k+s)!}{k!}.$ 
	Applying this above with $m=k-n,$ we have
	$$
	\lim\limits_{q\to 1} P(n,k)_q =(-1)^k\binom{n}{k}\sum\limits_{s=0}^{\lfloor\frac{n-k}{2}\rfloor} \frac{\left(\frac{k-n}{2}\right)_s\left(\frac{k-n+1}{2}\right)_s}{(k+1)_s}\cdot\frac{4^s}{s!},
	$$
	which is our statement since the summand vanishes for $s> \lfloor\frac{n-k}{2}\rfloor.$

	\subsection{Proof of Corollary~\ref{LegendreDistributionMatrices}}
	
	We prove this corollary by implementing the method of moments, as employed in previous work by the second two authors in \cite{hypergeometricDistribution}.
	By Theorem~\ref{EllipticCurveCase} and the fact that $p^{\frac{r}{2}}{_2F^{\ff}_1}(a)_{p^r}\in[-2,2],$ we have
	$$
	p^{\frac{r}{2}-rn^2}a_{\Leg,n}(a;p^r) = \phi_{p^r}(-1)\cdot p^{\frac{r}{2}}{_2F^{\ff}_1(a)_{p^r}} + O_{r,n}(p^{-\frac{r}{2}}).
	$$
	Therefore, if $m$ is a nonnegative integer, we have that
	\begin{align*}
	&\frac{1}{p^r}\sum\limits_{a\in\F_{p^r}\setminus\{0,1\}}\left(p^{\frac{r}{2}-rn^2}a_{\Leg,n}(a;p^r)\right)^m 
	=  \frac{1}{p^r}\sum\limits_{a\in\F_{p^r}\setminus\{0,1\}} (\phi_{p^r}(-1) p^{\frac{r}{2}}{_2F^{\ff}_1(a)_{p^r}} )^m \\ 
	&+ \frac{1}{p^r}\sum\limits_{k=1}^{m}\frac{1}{p^{\frac{rk}{2}}}\cdot \binom{m}{k}\frac{1}{p^r}\sum\limits_{a\in\F_{p^r}\setminus\{0,1\}} (\phi_{p^r}(-1) p^{\frac{r}{2}}{_2F^{\ff}_1}(a)_{p^r})^{m-k}+o_{m,r,n}(1) \\
	&= \frac{1}{p^r}\sum\limits_{a\in\F_{p^r}\setminus\{0,1\}} (\phi_{p^r}(-1) p^{\frac{r}{2}}{_2F^{\ff}_1(a)_{p^r}} )^m + o_{m,r,n}(1) \text{  as  }p\to\infty.
	\end{align*}
	By Theorem 1.1 of \cite{hypergeometricDistribution}, this implies that as $p\to\infty$ we have
	$$
	\frac{1}{p^r}\sum\limits_{a\in\F_{p^r}\setminus\{0,1\}}\left(p^{\frac{r}{2}-rn^2}a_{\Leg,n}(a;p^r)\right)^m = \begin{cases}
		o_{m,r,n}(1) & \text{   if   } m\text{ is odd } \\ 
		\frac{(2l)!}{l!(l+1)!} + o_{m,r,n}(1) & \text{  if   } m=2l\text{ is even }.
	\end{cases}
	$$
	The proof of Corollary 1.2 of \cite{hypergeometricDistribution} then implies the limiting distribution.

	\subsection{Proof of Theorem~\ref{K3Case}}
	Fix a prime power $q=p^r$ with $p\geq 3$ and $r\geq 1,$ and fix $a\in\F_q\setminus\{0,-1\}.$ If $A_{a}$ denotes the smooth affine surface given by
	$$
	s^2= xy(x+1)(y+1)(x+a y), s\neq 0,
	$$
	then $X:=A_a$ and $X_a$ differ by a connected union of rational curves (see \cite[\S 1]{ahlgrenOnoPenniston}).  In particular, we have
	$$
	[X_a] = [X] +24\mathbb{L}-6
	$$
	in the Grothendieck ring of $\F_q$-varieties, where $\mathbb{L}$ is the class of the affine line (cf. the term $|\psi^{-1}(U_a)|=(24q-6)$ at the end of the proof of \cite[Proposition 4.1]{ahlgrenOnoPenniston}). Therefore, by Theorem 1.1 of \cite{ahlgrenOnoPenniston}, the local zeta function of $X$ is given by 
	$$
	Z_{X}(t)=\frac{(1-qt)^5}{(1-t)^7(1-q^2t)(1-\gamma qt)(1-\gamma\alpha^2t)(1-\gamma\overline{\alpha}^2t)},
	$$
	where $\gamma=\phi_q(a+1)$ and $\alpha,\overline{\alpha}$ are the Frobenius eigenvalues for the Clausen elliptic curve $E_{\CL}\left(\frac{-a}{a+1}\right).$
	
	Therefore, by Proposition~\ref{prop:surface_zeta_formula}, we have that
	\begin{align*}
	\hat{Z}_X(t)  &=\prod\limits_{i,j\geq 1}\frac{(1-q^{1-j}t^i)^5}{(1-q^{-j}t^i)^7} \prod\limits_{i,j\geq 1} \frac{1}{(1-q^{2-j}t^{i})(1-\gamma q^{1-j}t^{i})(1-\gamma\alpha^2 q^{-j}t^{i})(1-\gamma\overline{\alpha}^2q^{-j}t^{i})} \\ 
	& = \prod\limits_{i\geq 1}(1-t^i)^5\cdot\prod\limits_{i,j\geq 1}\frac{1}{(1-q^{-j}t^i)^2}\cdot\prod\limits_{i\geq 1}\prod\limits_{b\in\{q^2,\gamma q,\gamma\alpha^2,\gamma\overline{\alpha}^2\}}\prod\limits_{j\geq 1} \frac{1}{1-bt^iq^{-j}}.
	\end{align*}
	By (\ref{Euler2}) and $\alpha\overline{\alpha}=q,$ we then have
	\begin{align*}
	&\hat{Z}_X(t) =(t;t)_\infty^5\cdot  \prod\limits_{i\geq 1}\left(\sum\limits_{m\geq 0}\frac{(-1)^mq^{\frac{m(m-1)}{2}}t^{mi}}{(q;q)_m}\right)^2\cdot\prod\limits_{i\geq 1}\prod\limits_{b\in\{q^2,\gamma q,\gamma{\alpha}^2,\gamma\overline{\alpha}^2\}}\sum\limits_{m\geq 0}\frac{(-1)^m q^{\frac{m(m-1)}{2}} b^m t^{im}}{(q;q)_m} \\
	&= (t;t)^5_\infty\prod\limits_{i\geq 1} \sum\limits_{m\geq 0} t^{im}\cdot\sum\limits_{m_1+\ldots+m_6=m}\frac{(-1)^m q^{\sum\frac{m_i(m_i-1)}{2}}\gamma^{m_4+m_5+m_6}\cdot q^{2m_3+m_4}\cdot \alpha^{2m_5}\overline{\alpha}^{2m_6}}{(q;q)_{m_1}\cdot\ldots\cdot (q;q)_{m_6}} \\ 
	&= (t;t)_\infty^5\sum\limits_{r\geq 0}t^r \cdot\sum\limits_{} (-1)^{\sum m_{u,v}} \frac{q^{\sum\frac{m_{u,v}(m_{u,v}-1)}{2}}q^{\sum 2m_{u,3}+m_{u,4}} \gamma^{\sum m_{u,4}+m_{u,5}+m_{u,6}} \cdot \alpha^{2\sum m_{u,5}}\overline{\alpha}^{2\sum m_{u,6}}}{\prod (q;q)_{m_{u,v}}},
	\end{align*}
	where the latter sum is over all possible combinations of nonnegative integers $m_{u,v}$ with $v=1,\ldots,6$ and positive integers $i_u$ such that $\sum i_{u}m_{u,v} = r.$ To simplify this expression, for each $v=1,\ldots,6,$  we denote by $\lambda_v$ the partition given by adding $i_u$ with multiplicity $m_{u,v}.$   
	
	Then, in the notation of theorem, the coefficient of $t^n$ in $\Zhat_X(t)$ is given by
	$$
	\sum\limits_{r=0}^nb_{n-r}\sum\limits_{\substack{\lambda_1,\ldots,\lambda_6 \\ |\lambda_1|+\ldots+|\lambda_6|=r}} (-1)^{l(\lambda_1)+\ldots+l(\lambda_6)} \frac{q^{\frac{\sum n(\lambda_i,j)(n(\lambda_i,j)-1)}{2}}q^{2l(\lambda_3)+l(\lambda_4)}}{\prod (q;q)_{n(\lambda_i,j)}}\gamma^{l(\lambda_4)+l(\lambda_5)+l(\lambda_6)} \alpha^{2l(\lambda_5)}\overline{\alpha}^{2l(\lambda_6)}.
	$$
	Dividing those partitions into $l(\lambda_5)-l(\lambda_6)=k$ for $0\leq k\leq n,$ using the definition\footnote{The set in \eqref{eq:matrix_point_k3} is indeed $C_n(X)$ as is defined in \eqref{eq:matrix_point_def}: a closed equation for $X$ is $s^2= xy(x+1)(y+1)(x+a y), sz=1$, so the corresponding matrix equation is $(A,B,C,D)\in C_{n,4}(\Fq), C^2=AB(A+I_n)(B+I_n)(A+aB), CD=I_n$, which is equivalent to \eqref{eq:matrix_point_k3} by setting $D=C^{-1}$.} of $\Zhat_X(t),$ the equation  $(\ref{Clausen3F2})$ and the relations $\alpha\overline{\alpha}=q$ and $\phi_{q^k}(a+1)=\phi_q(a+1)^k,$ we have that
	\begin{align*}
	N_{n,3}(a;q)&=S\left(n,0,\phi_q(a+1)\right)+\sum\limits_{k=1}^n S(n,k,\phi_q(a+1))_q\left((\alpha^k+\overline{\alpha}^k)^2-2\alpha^k\overline{\alpha}^k\right) \\
	&=S\left(n,0,\phi_q(a+1)\right)+\sum\limits_{k=1}^n S(n,k,\phi_q(a+1))_q\left( a_{\CL}\left(\frac{-a}{a+1}; q^k\right)^2-2q^k\right)\\
	&=S\left(n,0,\phi_q(a+1)\right)+\sum\limits_{k=1}^n S(n,k,\phi_q(a+1))_q\left(q^{2k}\cdot\phi_q(a+1)^k\cdot{_3F^{\ff}_2\left(-a\right)_{q^k}-q^k}\right),
	\end{align*}
	where
	\begin{align*}
	S(n,k,\gamma)_q := (-1)^n q^{\frac{n(n-1)}{2}}(q;q)_n\cdot  \sum\limits_{r=k}^n b_{n-r}\sum\limits_{\substack{|\lambda_1|+\ldots+|\lambda_6|=r \\ l(\lambda_5)-l(\lambda_6)=k}} & (-1)^{m(\lambda_1,\ldots,\lambda_6)} \gamma^{l(\lambda_4)+k} \\ 
	&\cdot\frac{q^{2l(\lambda_3)+l(\lambda_4)+2l(\lambda_6)} q^{\frac{\sum n(\lambda_i;j)(n(\lambda_i;j)-1)}{2}}}{\prod (q;q)_{n(\lambda_i,j)}}.
	\end{align*}
	The expression for $N_{n,3}(a;q)$ follows immediately. 
	
	To obtain the leading term for $Q(n,k,\gamma)_q,$ note that the leading term of $(q;q)_l$ is given by $(-1)^lq^{\frac{l(l+1)}{2}}.$ Therefore, the leading term of the summand for fixed $r,k$ and $\lambda_1,\ldots,\lambda_6$ is given by
	\begin{align*}
	&q^{\frac{n(n-1)}{2}}\cdot b_{n-r}\cdot q^{2l(\lambda_3)+l(\lambda_4)+2l(\lambda_5)}\gamma^{l(\lambda_4)}(-1)^{n-m(\lambda_1,\ldots,\lambda_6)}\cdot q^{\sum\frac{n(\lambda_i,j)(n(\lambda_i,j)-1)}{2}} \\
	&\cdot (-1)^n q^{\frac{n(n+1)}{2}} \cdot (-1)^{\sum n(\lambda_i,j)}q^{-\sum\frac{n(\lambda_i,j)(n(\lambda_i,j)+1)}{2}} \\
	&=b_{n-r}\gamma^{l(\lambda_4)}q^{n^2+2l(\lambda_3)+l(\lambda_4)+2l(\lambda_5)-m(\lambda_1,\ldots,\lambda_6)} \\
	&=b_{n-r}\gamma^{l(\lambda_4)}q^{n^2-l(\lambda_1)-l(\lambda_2)+l(\lambda_3)+l(\lambda_5)-l(\lambda_6)}.
	\end{align*}
	The exponent of $q$ is maximized when $l(\lambda_1)=l(\lambda_2)=l(\lambda_6)=0$ and $l(\lambda_3)+l(\lambda_5)$ is maximal. This implies that $l(\lambda_4)=0,$ $r=n$ and $\lambda_3$ and $\lambda_5$ are of the form $(1,\ldots,1).$ However, since $l(\lambda_5)-l(\lambda_6)=n,$ we thus have that $\lambda_5 = n\cdot 1,$ and therefore, the leading term of $Q(n,k,\gamma)_q$ is $q^{n^2+n}.$

	It remains to show the behavior of this polynomial as $q\to 1.$ 
	To this end, note that 
	$$
	\lim\limits_{q\to 1}(q;q)_{n-m(\lambda_1,\ldots,\lambda_6)}= 0
	$$
	if $m(\lambda_1,\ldots,\lambda_6)<n.$ Therefore, the only contributing partitions are those with $r=n$ and $l(\lambda_1)+\ldots+l(\lambda_6)=n.$ This implies that the partitions $\lambda_1,\ldots,\lambda_6$ are of the form $\lambda_i=(1,1,\ldots,1).$ Therefore,  since $b_0=1,$ we have
	\begin{align*}
		\lim\limits_{q\to 1} Q(n,k,\gamma)_q& = \sum\limits_{\substack{x_1+\ldots+x_6 = n \\ x_5-x_6 = k} }\frac{n!}{x_1!\ldots x_6!}\cdot\gamma^{x_4} \\ 
		&=\sum\limits_{x_1+\ldots+x_4+2x_6=n-k} \gamma^{x_4}\frac{n!}{x_1!\ldots x_4!x_6!(x_6+k)!} \\
		&=\sum\limits_{s=0}^{n-k}\left(\sum\limits_{x_1+x_2+x_3=s} \binom{n}{x_1,x_2,x_3,n-s}\right)\cdot\sum\limits_{x_4+2x_6=n-k-s}\binom{n-s}{x_4,x_6,x_6+k}\gamma^{x_4}
	\end{align*}
	The rest of the computation proceeds exactly as that of Theorem~\ref{EllipticCurveCase}.

	\subsection{Proof of Corollary~\ref{K3DistributionMatrices}}
		We prove this corollary by implementing the method of moments, as employed in previous work by the second two authors in \cite{hypergeometricDistribution}.
	By Theorem~\ref{K3Case} and the fact that $p^r{_3F^{\ff}_2}(a)_{p^r}\in[-3,3]$ which follows from (\ref{Clausen3F2}) we have
	$$
	p^{r-rn^2-rn}A_n(a;p^r) = p^r{_3F^{\ff}_2}(-a)_{p^r} + O_{r,n}(p^{-r}).
	$$
	Therefore, as in the proof of Corollary~\ref{K3DistributionMatrices}, and using Theorem 1.3 of \cite{hypergeometricDistribution}, for positive integers $m,$ we have that
	$$
	\frac{1}{p^r}\sum\limits_{a\in\F_{p^r}\setminus\{0,1\}} \left(p^{r-rn^2-rn}A_n(a;p^r)\right)^m = \begin{cases}
		o_{m,r,n}(1) & \text{   if  }m\text{ is odd} \\
		\sum\limits_{i=0}^m  (-1)^i \binom{m}{i}\frac{(2i)!}{i!(i+1)!} + o_{m,r,n}(1) & \text{ if }m\text{ is even.}
	\end{cases}
	$$
	The proof of Corollary 1.4 of \cite{hypergeometricDistribution} then implies the limiting distribution.

\bibliographystyle{amsplain}

\end{document}